\documentclass[reqno]{amsproc}
\usepackage[latin1]{inputenc}
\usepackage{amssymb}
\usepackage{hyperref}
\usepackage{amsmath}
\usepackage{latexsym}
\usepackage{cite}
\usepackage{amssymb}
\usepackage{amsfonts}
\usepackage{amscd}
\usepackage{xcolor}
\usepackage{amstext,amsmath,amssymb,amsfonts}
\newtheorem{theorem}{Theorem}

\newtheorem{corollary}[theorem]{Corollary}

\newtheorem{definition}[theorem]{Definition}
\newtheorem{example}[theorem]{Example}

\newtheorem{proposition}[theorem]{Proposition}
\newtheorem{remark}[theorem]{Remark}

\newcommand{\N}{\mathbb{N}}
\newcommand{\K}{\mathbb {K}}

\newcommand{\A}{\mathcal{A}}

\newcommand{\G}{\mathcal{G}}
\newcommand{\B}{\mathcal{B}}
\newcommand{\bb}{\mathfrak{B}}

\newcommand{\h}{\mathcal{H}}

\newcommand{\beq}{\begin{eqnarray}}
\newcommand{\eeq}{\end{eqnarray}}
\newcommand{\beqs}{\begin{eqnarray*}}
\newcommand{\eeqs}{\end{eqnarray*}}
\newcommand{\bpro}{\begin{pro}}
\newcommand{\epro}{\end{pro}}
\newcommand{\blem}{\begin{lem}}
\newcommand{\elem}{\end{lem}}
\newcommand{\bdfn}{\begin{dfn}}
\newcommand{\edfn}{\end{dfn}}
\newcommand{\bcor}{\begin{cor}}
\newcommand{\ecor}{\end{cor}}
\newcommand{\bthm}{\begin{thm}}
\newcommand{\ethm}{\end{thm}}
\newcommand{\bex}{\begin{ex}}
\newcommand{\eex}{\end{ex}}
\newcommand{\brmk}{\begin{rmk}}
\newcommand{\ermk}{\end{rmk}}
\newcommand{\bpr}{\begin{pr}}
\newcommand{\epr}{\end{pr}}
\newcommand{\benum}{\begin{enumerate}} 
\newcommand{\eenum}{\end{enumerate}}
\newcommand{\bitem}{\begin{itemize}}
\newcommand{\eitem}{\end{itemize}}

\newcommand{\cqfd}{\hfill{\square}}
\chardef\bslash=`\\
\numberwithin{equation}{section}
\numberwithin{table}{section}
\numberwithin{theorem}{section}
\DeclareMathOperator{\id}{id}
\DeclareMathOperator{\ad}{ad}

\title{Center-symmetric algebras  and bialgebras : relevant properties and   consequences}

\author{Mahouton Norbert Hounkonnou$^\ast$}
\address[$\ast$]{University of Abomey-Calavi,
International Chair in Mathematical Physics and Applications,
ICMPA-UNESCO Chair, 072 BP 50, Cotonou, Rep. of Benin}
\email{norbert.hounkonnou@cipma.uac.bj, with copy to hounkonnou@yahoo.fr}
\author{Mafoya Landry Dassoundo$^\dagger$}
\address[$\dagger$]{University of Abomey-Calavi,
International Chair in Mathematical Physics and Applications,
ICMPA-UNESCO Chair, 072 BP 50, Cotonou, Rep. of Benin}
\email{mafoya.dassoundo@cipma.uac.bj,  with copy to maflandas@gmail.com}
\begin{document}
\begin{abstract}
Lie admissible algebra structures, called
 center~-~symmetric algebras, are defined. 
Main properties and algebraic consequences are
derived and discussed. Bimodules are given
and used to build a center~-~symmetric
algebra on the direct sum of   underlying vector space 
and a finite dimensional  vector space.
Then,  the matched pair of center~-~symmetric
algebras is established and related to the
matched pair of sub-adjacent Lie
algebras. Besides, Manin triple of
center~-~symmetric algebras is defined
and linked with   their associated matched pairs.
Further, center~-~symmetric bialgebras of
center~-~symmetric algebras are investigated and discussed. Finally, a theorem  yielding the equivalence between Manin triple of center~-~symmetric algebras, matched pairs of Lie algebras and center-symmetric bialgebras is provided. 

\end{abstract}
\maketitle
\today
 


\tableofcontents
  \section{Introduction}
 Consider the algebra $(\A, \mu),$ i.e., a $\K$ vector space $\A$ endowed with a binary operation or law (bilinear homomorphism)  $\mu$  defined as:
 \beqs
 \mu: \A\times\A & \longrightarrow & \A \cr
 (x, y) &\longmapsto & \mu(x, y).
 \eeqs
Define the {\it associator  of the binary product} by   a trilinear homomorphism on $\A$  as follows \cite{rem.goze}:
 \beqs
 ass_\mu : \A\times\A\times\A &\longrightarrow & \A \cr
 (x, y, z) & \longmapsto & \mu(\mu(x,y), z)-\mu(x, \mu(y, z)).
 \eeqs
Let $\sigma \in \Sigma_3$ (symmetry group of degree $n$ $(n \in \N)$), acting on the associator 
 as: 
 \beqs
 \sigma (x_1, x_2, x_3) =\left(x_{\sigma^{-1}(1)}, x_{\sigma^{-1}(2)}, x_{\sigma^{-1}(3)}\right).
 \eeqs
 \begin{definition}\cite{santilli}
 The algebra $\A= (\A, \mu)$ is called Lie admissible if the commutator of $\mu,$ denoted by $[,]_{\mu},$ defines on $\A$ a Lie algebra structure, i.e., 
 $[x, y]_{\mu}=\mu(x, y)-\mu(y, x)$ (bilinear and skew-symmetric) and 
 $[[x, y]_{\mu}, z]_{\mu}+[[z, x]_{\mu}, y]_{\mu}+[[y, z]_{\mu}, x]_{\mu}=0$ (Jacobi identity).
 \end{definition}
 \begin{definition}\cite{rem.goze}
 The algebra $(\A, \mu)$ is called Lie~-~admissible if and only if $\mu$ satisfies:
 \beq
 \sum_{\sigma \in \Sigma_{3}} \left(-1\right)^{\varepsilon(\sigma)}ass_{\mu}\circ \sigma =0,
 \eeq
 where $\varepsilon$ is the signature of $\sigma.$
 \end{definition}
 \begin{definition}\cite{remm}
 Let $\mathrm{G}$ be a subgroup of $\Sigma_3.$ We say that the algebra law is $\mathrm{G}$~-~associative if 
 \beq
 \sum_{\sigma \in \mathrm{G}} \left(-1\right)^{\varepsilon(\sigma)}ass_{\mu}\circ \sigma =0.
 \eeq
 \end{definition}
 The subgroups of $\Sigma_3$ are well known. We have: $\mathrm{G}_1=\{\id \}$, $\mathrm{G}_2=\{\id, \tau_{12}\},$ $\mathrm{G}_3=\{\id, \tau_{23} \},$ $\mathrm{G}_4=\{ \id, \tau_{13} \}, $
 $\mathrm{G}_5=\{ \mathrm{A}_3 \}$ (Alternating group) and $\mathrm{G}_6=\Sigma_3.$  $\tau_{ij}$ is the transposition between $i$ and
 $j,$ i.e., explicitly:
  
 $\displaystyle \tau_{12}=\left(
 \begin{matrix} 
 1 & 2 & 3 \cr
 2 & 1 & 3
 \end{matrix}
 \right), 
 \tau_{13}=\left(
 \begin{matrix} 
 1 & 2 & 3 \cr
 3 & 2 & 1
 \end{matrix}
 \right), 
 \tau_{23}=\left(
 \begin{matrix} 
 1 & 2 & 3 \cr
 1 & 3 & 2
 \end{matrix}
 \right)$ and $\displaystyle \id=
  \left(
 \begin{matrix} 
 1 & 2 & 3 \cr
 1 & 2 & 3
 \end{matrix}
 \right).
 $
 
  We deduce the following types of Lie admissible algebras:
 \begin{enumerate}
 \item If $\mu$ is $\mathrm{G}_1$~-~associative, then $\mu$ is associative law. 
 \item If $\mu$ is $\mathrm{G}_2$~~-associative,
 then $\mu$ is a law of Vinberg algebra \cite{vinberg}. If $\A$ is finite~-~dimensional, then the associated Lie admissible algebra is provided with an affine structure.  
 \item If $\mu$ is $\mathrm{G}_3$~-~associative, then $\mu$ is a law of pre~-~Lie algebra (also called left~-~symmetric algebra).
 \item If $\mu$ is $\mathrm{G}_4$~-~associative, then $\mu$ satisfies
 \beq
 (xy)z-x(yz)=(zy)x-z(yx), \qquad \forall x, y, z \in \A. 
 \eeq
  We called the corresponding algebra {\it center-symmetric algebra}.
 \item If $\mu$ is $\mathrm{G}_5$ associative, then $\mu$ satisfies the generalized Jacobi condition i.e.
 \beq
 (xy)z+(yz)x+(zx)y=x(yz)+y(zx)+z(xy).
 \eeq
 Moreover, if the law is antisymmetric, then it is a law of Lie algebra.
 \item If $\mu$ is $\mathrm{G}_6$~-~associative, then $\mu$ is a Lie admissible Law.
 \end{enumerate}
 
 This work aims at  studying   $\mathrm{G}_4$~-~associative structures,  called center~-~symmetric algebras. Their  algebraic properties  are investigated.
  Related bimodule and  matched pairs  are given. Associated Manin triples built look like the Manin triple of Lie algebras \cite{Chari.Pressley}.
     Besides, from symmetry role of  matched pairs,  equivalent relations are  established in  the framework of center~-~symmetric bialgebras  making some bridges with the Lie~-~bialgebra construction by Drinfeld \cite{drinfeld}.
  
 Thoughout this work, we consider  $\A,$  a finite dimensional  vector space over the field $ \K $ of characteristic zero $(0)$ together with a bilinear product $\cdot$ defined as 
 $\displaystyle 
 \cdot : \A \times \A \rightarrow~\A$ such that $\; (~x~,~y~) \mapsto~x~\cdot~y.$
  
 \section{Basic properties: main definitions and algebraic consequences}
In this section, we give the definition of the center~-~symmetric algebra, provide their basic properties and deduce relevant algebraic consequences, similarly to known framework of  left~-~symmetric algebras \cite{bai1}.
 \begin{definition}
 $(\A, \cdot),$ (or simply $\A$), is said to be a center-symmetric algebra 
  if $\forall \, x, y, z \, \in \A,$  the associator of the bilinear product $\cdot,$ defined by $(~x,~y,~z~):= (~x\cdot y~)\cdot z -x\cdot~(~y\cdot z~),$ is symmetric in $x$ and $z$, i.e.,
 \begin{eqnarray}\label{associator}
 (~x,~y,~z~) = (~z,~y,~x~).
 \end{eqnarray}
 
 \end{definition} 
As matter of  notation simplification, we will denote $\displaystyle x \cdot y$ by $xy$ if not any confusion.
 \begin{remark}
 Any associative algebra is a center-symmetric algebra.
 \end{remark} 
\begin{proposition}
The bilinear product (commutator) 
$[\cdot , \cdot]:\, \A \times \A \longrightarrow  \A, \\ (x, y)  \longmapsto [x, y]=x\cdot y-y\cdot x$ gives a Lie bracket structure on $\A,$ known as the sub-adjacent Lie algebra $\G(\A):=(\A, [\cdot,\cdot])$ of  $(\A, \cdot).$
\end{proposition}
\textbf{ Proof:}
By definition of the commutator, $[\cdot, \cdot]$ is bilinear and skew symmetric. 
The Jacoby identity comes from a straightforward computation. $\cqfd$ 
%

 Thus, as in the case of lef-symmeric algebras, $(\A, \cdot)$ can be called the compatible center~-~symmetric algebra structure of the Lie algebra $\G(\A).$
 
 Considering the representations of the  left $L$ and right $R$ multiplication operations:
 \begin{eqnarray}
 L: \A & \longrightarrow & \mathfrak{gl}(\A)  \cr
  x  & \longmapsto & L_x:
  \begin{array}{ccc}
 \A &\longrightarrow & \A \cr 
  y & \longmapsto & x \cdot y, 
   \end{array}
\end{eqnarray}
\begin{eqnarray}
    R: \A & \longrightarrow & \mathfrak{gl}(\A)  \cr
     x  & \longmapsto & R_x:
     \begin{array}{ccc}
    \A &\longrightarrow & \A \cr 
     y & \longmapsto & y \cdot x,
      \end{array}
 \end{eqnarray}
 we infer the adjoint representation  $\ad: = L-R$ of the sub-adjacent Lie algebra $\G(\A)$ of a center~-~symmetric algebra $\A$  as follows:
  \begin{eqnarray}
    \ad: \A & \longrightarrow & \mathfrak{gl}(\A)  \cr
      x  & \longmapsto & \ad_x:
      \begin{array}{ccc}
     \A &\longrightarrow & \A \cr 
      y & \longmapsto & [x, y],
       \end{array}
    \end{eqnarray}
such that  
  $\displaystyle \forall \; x, y \in \A, \ad_x(y):=(L_x-R_x)(y).$ 
 \begin{proposition}\label{one}
 Let $(\A, \cdot)$ be a center~-~symmetric algebra, and $L$, (resp. $R$), be the linear representation of the left, (resp. right), multiplication operator. Then,
 \begin{enumerate}
\item For all $\displaystyle x, y \in \A$ we have: 
$\displaystyle [L_x, R_y]=[L_y, R_x]$ and 
$\displaystyle L_{x\cdot y}-L_xL_y=R_xR_y-R_{y\cdot x}.$
\item $\ad= L-R$ is a linear representation of the sub-adjacent Lie algebra $\G(\A)$ of $(\A, \cdot),$ i.e., \\
$\displaystyle \ad_{[x, y]}=[\ad_x, \ad_y], \; \forall \;  x,  y \, \in \A.$ 
 \end{enumerate}
 \end{proposition}
\textbf{ Proof:} It is immediate from the definitions of considered operators. $\cqfd$
 \section{Bimodules and matched pairs}
  \begin{definition}\label{bimodule}
Let $\A$ be a center~-~symmetric algebra, $\displaystyle V$ be a vector space. Suppose  \\ $\displaystyle l,r : \A \rightarrow \mathfrak{gl}(V)$ be two linear maps satisfying: 
$\mbox{For all } x, y \in \A,$
\begin{eqnarray}\label{eqbimodule1}
\left[l_x, r_y\right]= \left[l_y, r_x\right]
\end{eqnarray}
\begin{eqnarray}\label{eqbimodule2}
l_{xy}-l_xl_y=r_xr_y-r_{yx}.
\end{eqnarray}
Then, $\displaystyle(l, r, V)$ (or simply $\displaystyle(l,r)$) is called bimodule of the center~-~symmetric algebra~$\A.$
 \end{definition}
In this case, the following statement can be proved by a direct computation. 
\begin{proposition}
Let $(\A, \cdot)$ be a center-symmetric algebra and $V$ be a vector space over $\K.$
Consider two linear maps,
 $\displaystyle l, r : \A \rightarrow \mathfrak{gl}(V).$ Then, $(l, r, V)$ is a bimodule of $\A$ if and only if, the semi-direct sum  $\A\oplus V $ of vector spaces is turned into a center-symmetric algebra  by defining the multiplication in $\A\oplus V$ by
$ \displaystyle
(x_1+v_1)\ast (x_2+v_2) = x_1\cdot x_2+(l_{x_1}v_2+r_{x_2}v_1), \; \forall x_1, x_2 \in \A, \; v_1, v_2 \in V.
$
We denote it by $\displaystyle \A \ltimes_{l, r} V $ or simply $\A \ltimes V.$
\end{proposition} 
Furthermore, we derive the next result.
\begin{proposition}\label{propcentlie}
Let  $\A$ be a center-symmetric algebra and $V$ be a vector space over $\K$.
Consider two linear maps,
 $\displaystyle l, r: \A \rightarrow \mathfrak{gl}(V),$ 
 such that 
$(l,r, V)$ is a bimodule of $\A.$ Then, the map:
$ \displaystyle 
 l-r: \A  \longrightarrow  \mathfrak{gl}(V) \;  
  x   \longmapsto  l_x-r_x,
  $
is a linear representation of the sub-adjacent Lie algebra of $\A.$
\end{proposition}
\begin{example}
According to the Proposition~\ref{one}, one can deduce that $(L, R, V)$ is a bimodule of the center-symmetric algebra $\A$, where $L$ and $R$ are the left and right multiplication operator representations, respectively.
\end{example}    
\begin{definition}\label{dfnmatched}\cite{majid}
Let $\displaystyle \mathcal{G}$ and $\displaystyle \mathcal{H}$ be two Lie algebras and let  $\displaystyle \mu: \mathcal{H} \rightarrow \mathfrak{gl}(\mathcal{G})$ and $\displaystyle \rho: \mathcal{G}\rightarrow \mathfrak{gl}( \mathcal{H})$ be two Lie algebra representations satisfying:
$\displaystyle \mbox{For all } x, y \in \mathcal{G}, a, b \in~\mathcal{H},$ 
\begin{eqnarray}\label{eqt1}
\rho(x)\left[a, b\right]-\left[\rho(x)a, b\right]-\left[a, \rho(x)b\right]+\rho(\mu(a)x)b-\rho(\mu(b)x)a=0,
\end{eqnarray}
\begin{eqnarray}\label{eqt2}
\mu(a)\left[x, y\right]-\left[\mu(a)x, y\right]-\left[x, \mu(a)y\right]+\mu(\rho(x)a)y-\mu(\rho(y)a)x=0.
\end{eqnarray}
Then, $\displaystyle(\mathcal{G}, \mathcal{H}, \rho, \mu)$ is called a matched pair of the Lie algebras
$\displaystyle \mathcal{G}$ and $\displaystyle \mathcal{H},$ 
denoted by $\displaystyle \mathcal{H}\bowtie_{\mu}^{\rho}\mathcal{G}.$
In this case,  $\displaystyle (\mathcal{G}\oplus \mathcal{H}, \ast)$ defines a Lie algebra with respect to the product $\ast$ satisfying: \\ $\displaystyle (x+a)\ast (y+b)= [x, y]+\mu(a)y-\mu(b)x+[a, b]+\rho(x)b-\rho(y)a$.
\end{definition}
 \begin{theorem}\label{theoo}
 Let $(\A, \cdot)$ and $(\B, \circ)$ be two center~-~symmetric algebras. 
Suppose that $(l_{\A}, r_{\A}, \B)$ and $(l_{\B}, r_{\B}, \A)$
 are bimodules of $\A$ and $\B$, respectively, obeying the relations:
\begin{eqnarray}\label{eqq1}
  -r_{\A}(x)(a \circ b)+r_{\A}(l_{\B}(b)x)a+ a\circ (r_{\A}(x)b)+l_{\A}(r_{\B}(b)x)a+ (l_{\A}(x)b)\circ a -l_{\A}(x)(b\circ a)=0,
\end{eqnarray}
\begin{eqnarray}\label{eqq2}
  -r_{\B}(a)(x \cdot y)+r_{\B}(l_{\A}(y)a)x+ x \cdot(r_{\B}(a)y)+l_{\B}(r_{\A}(y)a)x+   (l_{\B}(a)y)\cdot x - l_{\B}(a)(y \cdot x) = 0, 
\end{eqnarray}
\begin{eqnarray}\label{eqq3}
  a \circ (l_{\A}(x)b)+(r_{\A}(x)b)\circ a -(r_{\A}(x)a)\circ b-l_{\A}(l_{\B}(a)x)b+  r_{\A}(r_{\B}(b)x)a+l_{\A}(l_{\B}(b)x)a+ \cr -b\circ (l_{\A}(x)a)  -r_{\A}(r_{\B}(a)x)b=0, 
 \end{eqnarray}
\begin{eqnarray}\label{eqq4}
 x\cdot (l_{\B}(a)y)+(r_{\B}(a)y)\cdot x - (r_{\B}(a)x)\cdot y-l_{\B}(l_{\A}(x)a)y+  r_{\B}(r_{\A}(y)a)x+l_{\B}(l_{\A}(y)a)x+\cr 
 -y\cdot (l_{\B}(a)x)-r_{B}(r_{\A}(x)a)y=0,
 \end{eqnarray}
 for any $x, y \in A$ and $a, b \in \B.$ Then, there is a center~-~symmetric algebra structure on $\A \oplus \B$ 
 given by: 
 $\displaystyle (x+a)\ast (y+b)= (x \cdot y + l_{\B}(a)y+r_{\B}(b)x)+ (a \circ b + l_{\A}(x)b+r_{\A}(y)a).$
 We denote this center~-~symmetric algebra by $\displaystyle \A \bowtie_{l_{\B}, r_{\B}}^{l_{\A}, r_{\A}} \B,$ or
 simply by $\displaystyle \A \bowtie \B.$ Then $(\A, \B, l_{\A}, r_{\A}, l_{\B}, r_{\B})$ satisfying
the above conditions is called  matched pair of the center~-~symmetric algebras $\A$ and $\B.$
 \end{theorem}
\textbf{Proof:} Consider $x, y \in \mathcal{A}$ and $a, b \in \mathcal{B}$. 
\beqs \mbox{We have \,   }
(x+a) \ast (y+b) &=& (x y+ l_{\mathcal{B}}(a)y+r_{\mathcal{B}}(b)x)+ (a \circ b + l_{\mathcal{A}}(x)b+r_{\mathcal{A}}(y)a), 
\eeqs
and the associator takes the form:
\beqs
(x+a, y+b, z+c) &=& (x, y, z) + (a, b, c) + \{r_{\mathcal{B}}(c)(x \cdot y) + 
l_{\mathcal{A}}(x \cdot y)c -x \cdot (r_{\mathcal{B}}(c)y) 
\cr
&&
-l_{\mathcal{A}}(x)(l_{\mathcal{A}}(y)c)
-r_{\mathcal{B}}(l_{\mathcal{A}}(y)c)x \} + 
\{r_{\mathcal{B}}(c)(r_{\mathcal{B}}(b)x) 
+ l_{\mathcal{A}}(r_{\mathcal{B}}(b)x)c 
\cr  && 
-r_{\mathcal{B}}(b \circ c)x + (l_{\mathcal{A}}(x)b)\circ c-
l_{\mathcal{A}}(x) ( b \circ c) \} 
+ \{(r_{\mathcal{B}}(b)x) \cdot z 
 l_{\mathcal{B}}(l_{\mathcal{A}}(x)b)z
 \cr &&
 + r_{\mathcal{A}}(z)(l_{\mathcal{A}}(x)b)-
x \cdot (l_{\mathcal{B}}(b)z)- 
 r_{\mathcal{B}}(r_{\mathcal{A}}(z)b)x-l_{\mathcal{A}}(x)(r_{\mathcal{A}}(z)b)\} +
\cr &&
 \{(l_{\mathcal{B}}(a)y)\cdot z
 + 
l_{\mathcal{B}}(r_{\mathcal{A}}(y)a)z 
+r_{\mathcal{A}}(z)(r_{\mathcal{A}}(y)a)
- l_{\mathcal{B}}(a)(y \cdot z)- r_{\mathcal{A}}(y \cdot z )a\}
\cr &&
+
\{r_{\mathcal{B}}(c)(l_{\mathcal{B}}(a)y) 
+(r_{\mathcal{A}}(y)a) \circ c +
 l_{\mathcal{A}}(l_{\mathcal{B}}(a)y)c-
+l_{\mathcal{B}}(a)(r_{\mathcal{B}}(c)y)
\cr &&
- a \circ (l_{\mathcal{A}}(y)c) 
-r_{\mathcal{A}}(r_{\mathcal{B}}(c)y)a\} + \{l_{\mathcal{B}}(a \circ b)z + r_{\mathcal{A}}(z)(a \circ b)-
\cr &&
+ l_{\mathcal{B}}(a)(l_{\mathcal{B}}(b)z) 
- a \circ (r_{\mathcal{A}}(z)b)- r_{\mathcal{A}}(l_{\mathcal{B}}(b)z)a\},
\eeqs
which can also be re-expressed as:
\beq
(x+a, y+b, z+c)&= &(x, y, z)+(x, y, c)+(x, b, z)+ (x, b, c) \cr
&+& 
 (a, y, z)+(a, y, c)+(a, b, z) +(a, b, c).
\eeq
Similarly, 
\beq
(z+c, y+c, x+a) &=& (z, y, x)+(z, y, a)+(z, b, x)+(z, b, a) \cr &+&(c, y, x)+(c, b, a)+(c, y, a)  +(c, b, x).
\eeq
Using the fact that $(l_{\mathcal{A}}, r_{\mathcal{A}})$ is a bimodule of $\mathcal{A}$ and $(l_{\mathcal{B}}, r_{\mathcal{B}})$ is 
a bimodule of 
$\mathcal{B},$ one arrives at the following result:
\beqs
(x+a, y+b, z+c) = (z+c, y+b, x+a)
 &\Longleftrightarrow&  
\left\lbrace 
\begin{array}{ccc} 
(x, y, z) &=& (z, y, x) \\
(x, y, c) &=& (c, y, x) \\
(x, b, z) &=& (z, b, x) \\
(x, b, c) &=& (c, b, x) \\
(a, y, z) &=& (z, y, a) \\
(a, y, c) &=& (c, y, a) \\
(a, b, z) &=& (z, b, a) \\
(a, b, c) &=& (c, b, a)\\
 \end{array}
\right.  \\
&\Longleftrightarrow& 
\left\lbrace
\begin{array}{ccc}
(l_{\mathcal{A}}, r_{\mathcal{A}}, \mathcal{B})  &,&
(l_{\mathcal{B}}, r_{\mathcal{B}}, \mathcal{A}) \\
(x, y, c) &=& (c, y, x)  \ \eqref{eqq2}\\
 (x, b, z) &=& (z, b, x)  \ \eqref{eqq4}\\
(x, b, c) &=& (c, b, x)  \ \eqref{eqq1}\\
 (a,y,c) &=& (c,y,a).\ \ \eqref{eqq3}\\				
\end{array}
	\right.    				 
\eeqs
This last relation  ends the proof.
$\cqfd$

Moreover, every center~-~symmetric algebra which 
is a direct sum of the underlying spaces of two center~-~symmetric sub~-~algebras can be obtained in the 
above way.  
\begin{corollary}\label{corr}
Let $(\A, \B, l_{\A}, r_{\A}, l_{\B}, r_{\B})$ be 
a matched pair of center~-~symmetric algebras. Then, 
$(~\mathcal{G}(\A), \mathcal{G}(\B), l_{ \A}-r_{ \A} , l_{\B}-r_{\B} ~)$ is a matched pair of sub-adjacent Lie algebras $\G(\A)$ and~$\G(\B).$
\end{corollary}
\textbf{Proof:} 
 By using the Proposition~\ref{propcentlie} and the bimodules
$\displaystyle (l_{\A}, r_{\A}, \B)$ and  $(l_{\B}, r_{\B}, \A),$ 
 we have: 
$\ad_{\A}:=l_{\A}-r_{\A}$ and $\ad_{\B}:=l_{\B}-r_{\B}$ are   the linear representations of the sub-adjacent Lie algebras $\G(\A)$ and $\G(\B)$ of the center-symmetric algebras $\A$ and $\B,$ respectively.
Then, the statement that
$\displaystyle\mathcal{G}{(\A)}\bowtie_{\ad_{\B}}^{\ad_{\A}}\mathcal{G}{(\B)}$ is a matched pair of the Lie algebras $\G(\A)$ and $\G(\B)$ follows from Theorem~\ref{theoo}. 
Hence, $(~\mathcal{G}(\A), \mathcal{G}(\B), \ad_{\A} , \ad_{\B})$ is a matched pair of sub-adjacent Lie algebras $\G(\A)$ and $\G(\B)$.$\cqfd$
\begin{definition}\label{dfnn}
Let $(l, r, V)$ be a bimodule of a center~-~symmetric algebra $\A,$ where $V$ is a finite dimensional vector 
space. The dual maps $l^{*}, r^{*}$  of the linear maps $l, r,$ are defined, respectively, as: 
$\displaystyle l^{*}, r^{*}: \A \rightarrow \mathfrak{gl}(V^{*})$
 such that:
\beq\label{dual1}
 l^*: \A & \longrightarrow & \mathfrak{gl}(V^*)  \cr
  x  & \longmapsto & l^*_x:
      \begin{array}{llll}
 V^* &\longrightarrow & V^* \\ 
  u^* & \longmapsto & l^*_x u^*: 
      \begin{array}{llll}
V  &\longrightarrow&  \K \cr
v  &\longmapsto& \left< l^{*}_xu^{*}, v \right>
 := \left<  u^{*}, l_x v\right>, 
      \end{array}
     \end{array}
 \eeq
\beq\label{dual2}
 r^*: \A & \longrightarrow & \mathfrak{gl}(V^*)  \cr
  x  & \longmapsto & l^*_x:
      \begin{array}{llll}
 V^* &\longrightarrow & V^* \\ 
  u^* & \longmapsto & r^*_x u^*: 
      \begin{array}{lllll}
V  &\longrightarrow&  \K \cr
v  &\longmapsto& \left< r^{*}_xu^{*}, v \right>
 := \left<  u^{*}, r_x v\right>, 
      \end{array}
     \end{array}
\eeq
for all $x \in \A, u^{*} \in V^{*}, v \in V.$
\end{definition}  
\begin{proposition}
Let  $\A$ be a center~-~symmetric algebra and   
$\displaystyle l, r: \A \rightarrow \mathfrak{gl}(V)$ be two linear maps, where $V$ is a finite dimensional vector space.
The following conditions are equivalent:
\begin{enumerate}
\item $(l, r, V)$ is a bimodule of $\A$.
\item $(r^{*}, l^{*}, V^{*})$ is a bimodule of
 $\A .$  
\end{enumerate}
\end{proposition}
{\textbf{Proof}:} It  stems from the Definition~\ref{dfnn}.
\begin{theorem}\label{ttheo}
Let $(\A, \cdot )$ be a center~-~symmetric algebra. Suppose that there exists a center~-~symmetric algebra structure $"\circ "$ on its dual space $\A^{*}.$ Then,  $(\A, \A^{*}, R_{\cdot}^*, L_{\cdot }^*, R_{\circ}^*, L_{\circ}^*)$ is a matched pair of center~-~symmetric algebras $\A$ and $\A^{*}$
if and only if $(~\G(\A), \G(\A^{*}), -\ad^*_{\cdot}, -\ad^{*}_{\circ}~)$ is a matched pair of Lie algebras $\G(\A)$ and $\G(\A^{*}).$
\end{theorem}
{\textbf{Proof}:}
By considering the Theorem~\ref{theoo}, setting   
$\displaystyle l_{\A}:=R^{*}_{\cdot}, r_{\A}:=L^*_{\cdot}, l_{\B}:=R^*_{\circ},  
r_{\B}:=L^*_{\circ},$ and exploiting the
 Definition~{\ref{dfnmatched}}
  with $\G:=\G(\A),
   \h:=\G(\A^*),\rho:=R_{\cdot}^*-L_{\cdot}^*, \mu:=R_{\circ}^*-L_{\circ}^*,$ and  the relations
\eqref{dual1} and \eqref{dual2}, we get the equivalences. $\cqfd$
\begin{proposition}
Let $\G$ be a  Lie algebra.  Suppose  $\displaystyle \rho : \G \rightarrow
 \mathfrak{gl}(V)$ and $\displaystyle   
 \mu: \G \rightarrow \mathfrak{gl}(W)$ 
 be two linear representations of $\G,$ where $V$ and $W$ are  two vector spaces.
Then, the linear map 
$\displaystyle \rho\otimes 1+1\otimes \mu : \G \rightarrow \mathfrak{gl}(V \otimes W)$ given by  
$(\rho\otimes 1+1\otimes \mu) (v,w):=\rho(x)v\otimes w+v\otimes \mu(x)v$ 
is also a representation of $\G.$
\end{proposition}
{\textbf {Proof}:} It comes from a
straightforward computation. $\cqfd$
 \begin{theorem}\label{theo}
 Let $\A$ be a center~-~symmetric algebra with the product given by the
 linear map $\beta^{*}: \A\otimes\A \rightarrow \A.$
 Suppose there is a center~-~symmetric algebra structure $"\circ"$ on the dual space $\A^{*}$ provided by a linear map 
 $\alpha^{*}: \A^{*}\otimes \A^{*} \rightarrow \A^{*}.$ Then, $(\mathcal{G}(\A), \mathcal{G}(\A^*), -\ad_{\cdot}^*, -\ad^{*}_{\circ })$ is a matched pair of Lie algebras $\G(\A)$ and $\G(\A^*)$ if and only if $\alpha : \A \rightarrow \A \otimes \A$
 is a $1$~cocycle of $\mathcal{G}(\A)$ associated to $(-\ad_{\cdot})\otimes1+1\otimes (-\ad_{\cdot })$ and $\beta : \A^{*} \rightarrow \A^{*}\otimes \A^{*}$ is a 1-cocycle of  $\mathcal{G}(\A^{*})$ associated to
 $(-\ad_{\circ})\otimes 1 + 1\otimes (-\ad_{\circ}).$
 \end{theorem}
{\textbf{Proof}:} See Appendix.
   \section{Manin triple and center-symmetric bialgebras}
 In this section, similarly to the notion of Manin triple of Lie
 algebras inroduced in \cite{Chari.Pressley}, we first give the
 definition of Manin triple of a center~-~symmetric
algebra and investigate its associated  bialgebra structure. 
 Then, 
   we provide the  
basic definition and properties of center~-~symmetric bialgebras.
\begin{definition}
 A Manin triple of  center~-~symmetric algebras is a triple $(\A, \A^{+}, \A^{-})$ 
 together with a   nondegenerate symmetric bilinear form $\bb~(~;~)$ on $\A$ which is invariant, i.e., $\forall x, y, z \in \A,$
 $\bb(x \ast y, z) = \bb(x, y\ast z),$ satisfying:
 \begin{enumerate}
 \item $\A=\A^{+}\oplus \A^{-}$ as 
 $\K$~-~ vector space;
 \item $\A^{+}$ and $\A^{-}$ are center~-~symmetric subalgebras of $\A;$ 
 \item  $\A^{+}$ and $\A^{-}$ are isotropic with respect  to 
 $\bb(;),$ i.e., $\bb(\A^{+};\A^{+})=0= \bb(\A^{-};\A^{-}).$
 \end{enumerate} 
\end{definition}
Two Manin triples $(\A_1, \A_1^{+}, \A_1^{-}, \bb_1)$ and $(\A_2, \A_2^{+}, \A_2^{-}, \bb_2)$ of center~-~symmetric algebras $\A_1$ and $\A_2$ are homomorphic (isomorphic) if there is a homomorphism (isomorphism) $\displaystyle \varphi: \A_1 \rightarrow \A_2$ such that:
$ \displaystyle 
\varphi(\A_1^{+})\subset \A_2^{+},\; \varphi(\A_1^{-})\subset \A_2^{-}, \;
\bb_1(x, y)=\varphi^*\bb_2(\varphi(x), \varphi(y))=\bb_2(\varphi(x), \varphi(y)).
$
In particular, if $(\A, \cdot)$ is a center~-~symmetric algebra, 
and if  there exists a center~-~symmetric algebra structure on its
dual space $\A^*$ denoted $(\A^*, \circ),$ then there is a
center~-~symmetric algebra structure on the direct sum of the
underlying vector space of $\A$ and  $\A^*$ 
(see Theorem~\ref{theoo} ) such that $(\A\oplus\A^*, \A, \A^*)$ 
is the associated Manin triple
  with the invariant bilinear symmetric form given by
  $ \displaystyle 
\bb_{\A}(x+a^*, y+b^*)=<x, b^*>+<y, a^*>, \; \forall x, y \in \A; a^*, b^* \in \A^*,
$
 called the standard Manin triple of the center~-~symmetric algebra $\A$.
\begin{theorem}\label{thheo}
Let $(\A, \cdot )$ and $(\A^*, \circ )$ be two center~-~symmetric algebras. Then, the sixtuple
 $(\A, \A^*, R_{\cdot}^*, L_{\cdot}^*; R_{\circ}^*, L_{\circ}^* )$ is a matched pair of center~-~symmetric algebras $\A$ and $\A^{*}$ if and only if \\  $(\A \oplus \A^*, \A, \A^* )$ is their standard Manin triple.
\end{theorem}
{\textbf{Proof}:}
By considering that $(\A, \A^*, R_{\cdot}^*, L_{\cdot}^*; R_{\circ}^*, L_{\circ}^* )$ is a matched pair of center~-~symmetric algebras, it follows that the bilinear product $\ast$ defined in the Theorem~\ref{theoo} 
is  center~-~symmetric  on the direct sum of underlying vectors spaces,  $\A \oplus \A^*.$ 
Computing and comparing the  relations, we get: $\displaystyle \bb_{\A}\left((x+a)\ast(y+b), (z+c)\right)=\displaystyle \bb_{\A}\left((x+a), (y+b)\ast(z+c)\right)$  $\forall x, y, z \in \A; a, b, c \in \A^{*},$
which expresses the invariance of the standard bilinear form on $\A \oplus \A^*.$ 
Therefore, $(A\oplus \A^*, \A, \A^*)$ is the standard Manin triple of the center~-~symmetric algebras $\A$ and  $\A^*.$
$\cqfd$
\begin{definition}
Let $\A$ be a vector space. A center~-~symmetric bialgebra structure on $\A$ is a pair of linear maps $(\alpha, \beta)$ such that $\alpha: \A \rightarrow \A \otimes \A,$ $\beta: \A^* \rightarrow \A^*\otimes\A^*$ satisfying:
\begin{enumerate}
\item $\alpha^*: \A^*\otimes\A^* \rightarrow \A^*$ is a center~-~symmetric algebra structure on $\A^*,$
\item $\beta^*: \A\otimes\A \rightarrow \A$ is a center~-~symmetric algebra structure on $\A,$
\item $\alpha $
 is a $1$~cocycle of $\mathcal{G}(\A)$ associated to $(-\ad_{\cdot})\otimes1+1\otimes (-\ad_{\cdot }),$ 
 \item $\beta$ is 
  1-cocycle of  $\mathcal{G}(\A^{*})$ associated to
 $(-\ad_{\circ})\otimes 1 + 1\otimes (-\ad_{\circ}).$
\end{enumerate}
We also denote this center~-~symmetric bialgebra by $(\A, \A^*, \alpha, \beta)$ or simply $(\A, \A^*).$
\end{definition}
\begin{proposition}
Let $(\A, \cdot)$ be a center~-~symmetric algebra and $(\A^*, \circ)$ be a center~-~symmetric algebra structure on its dual space $\A^*.$ Then the following conditions are equivalent:
\begin{enumerate}
\item $(\A\oplus\A^*, \A, \A^*)$ is the standard Manin triple of considered center~-~symmetric algebras; 

\item $(\G(\A), \G(\A^*),-\ad_{\cdot}^*, -\ad_{\circ}^*)$ is a  matched pair of sub-adjacent Lie algebras;
\item $(\A, \A^*, R_{\cdot}^*, L_{\cdot}^*, R_{\circ}^*, L_{\circ}^*)$ is a matched pair of center~-~symmetric algebras;
\item $(\A, \A^*)$ is a center~-symmetric bialgebra.
\end{enumerate}
\end{proposition}
{\textbf{Proof}:}
From Theorem~\ref{ttheo},
$\displaystyle (2) \Longleftrightarrow (3),$ while from  Theorem~\ref{theo},  
$\displaystyle (2) \Longleftrightarrow (4).$
 Theorem~\ref{thheo} shows that 
$(1) \Longleftrightarrow (3).$
$\cqfd$
\section{Concluding remarks}
In this work, we have defined Lie admissible algebra structures, called
 center~-~symmetric algebras for which  
main properties and algebraic consequences have been
derived and discussed. Bimodules have been given
and used to build a center~-~symmetric
algebra on the direct sum of a
 center~-~symmetric algebra and a vector space.
Then, we have established the matched pair of center~-~symmetric
algebras, which has been  related to the
matched pair of sub~-~adjacent Lie
algebras. Besides, we have defined the Manin triple of
center~-~symmetric algebras 
and linked it with   their associated matched pairs.
Further, we have investigated   and discussed center~-~symmetric bialgebras of
center~-~symmetric algebras. Finally, we have provided a theorem  yielding the equivalence between Manin triple of center~-~symmetric algebras, matched pairs of Lie algebras and  center~-~symmetric algebras, and center~-~symmetric bialgebra.

\section*{Appendix}
{\textbf{Proof of the Theorem~\ref{theo}}} \\
Let $\{e_1, e_2, \cdots, e_n\}$ be a basis of $\A$ and  $\{e_1^{*}, e_2^{*}, \cdots, e_n^{*}\}$ 
 its dual basis.
Consider
$\displaystyle e_i\cdot e_j = \sum_{k=1}^{n}c_{ij}^{k}e_k$ and $\displaystyle e_i^{*}\circ e_j^{*} 
= \sum_{k=1}^{n}f_{ij}^{k}e_k^{*},$ where $\displaystyle c_{ij}^k, f_{ij}^k \in \K$ are  structure constants associated to $\cdot$ and $\circ,$ respectively. Then, it   follows that:\\
$\displaystyle \alpha(e_k) = \sum_{i,j=1}^{n}f_{ij}^{k}e_i\otimes e_j,$ 
$\displaystyle \beta(e_k^{*})= \sum_{i,j=1}^{n}c_{ij}^{k}e_i^{*}\otimes e_j^{*},$
and
\beq\label{cocycle1}
\alpha([e_i, e_j]) = \sum_{m,l=1}^{n}\sum_{k=1}^{n}\left\lbrace (c_{ij}^{k}-c_{ji}^{k})f_{ml}^{k} \right\rbrace e_m \otimes e_l,
\eeq
\beq\label{dualcocycle1}
\beta([e_i^{*}, e_j^{*}]) = \sum_{m,l=1}^{n}\sum_{k=1}^{n}\left\lbrace (f_{ij}^{k}-f_{ji}^{k})c_{ml}^{k} \right\rbrace e_m^{*}\otimes e_l^{*},
\eeq
and we get:
\beq\label{cocycle2}
\{(-\ad_{\cdot})(e_i)\otimes 1+1\otimes (-\ad_{\cdot})(e_i)\}\alpha(e_j)-
\{(-\ad_{\cdot})(e_j)\otimes1+1\otimes (-\ad_{\cdot})(e_j)\}\alpha(e_i) =\cr
\sum_{m,l=1}^{n}\sum_{k=1}^{n}\left\lbrace -f_{kl}^{j}(c_{ik}^{m}-c_{ki}^{m})+f_{kl}^{i}(c_{jk}^{m}-c_{kj}^{m})-f_{mk}^{j}(c_{ik}^{l}-c_{ki}^{l})
+f_{mk}^{i}(c_{jk}^{l}-c_{kj}^{l})\right\rbrace e_m\otimes e_l
\eeq
Taking into account  the fact that $\alpha$ is a $1$-cocycle of $\G(\A)$ associated to  $(-\ad_{\cdot})\otimes 1 + 1 \otimes(-\ad_{\cdot})$, and using the relations
\eqref{cocycle1} and \eqref{cocycle2} yield:
\beq\label{A}
\sum_{k=1}^{n} (c_{ij}^{k}-c_{ji}^{k})f_{ml}^{k}= \sum_{k=1}^{n} \left\lbrace f_{kl}^{i}(c_{jk}^{m}-c_{kj}^{m}) -f_{kl}^{j}(c_{ik}^{m}-c_{ki}^{m})+ 
 f_{mk}^{i}(c_{jk}^{l}-c_{kj}^{l})- 
 f_{mk}^{j}(c_{ik}^{l}-c_{ki}^{l} )\right\rbrace.
\eeq
Besides, we obtain:
\beq\label{dualcocycle2}
\{(-\ad_{\circ})(e_i^{*})\otimes1+1\otimes (-\ad_{\circ})(e_i^{*})\}\beta(e_j^{*})-
\{(-\ad_{\circ})(e_j^{*})\otimes 1+ 1\otimes (-\ad_{\circ})(e_j^{*})\}\beta(e_i^{*})= \cr
\sum_{m,l=1}^{n}\sum_{k=1}^{n}\left\lbrace-c_{kl}^{j}(f_{ik}^{m}-f_{ki}^{m})+c_{kl}^{i}(f_{jk}^{m}-f_{kj}^{m})- 
c_{mk}^{j}(f_{ik}^{l}-f_{ki}^{l})+c_{mk}^{i}(f_{jk}^{l}-f_{kj}^{l})\right\rbrace(e_m^{*}\otimes e_{l}^{*}).
\eeq
As $\beta$ is the $1$-cocycle issued from   $(-\ad_{\circ})\otimes1+1\otimes(-\ad_{\circ})$ and using the relations
\eqref{dualcocycle1} and \eqref{dualcocycle2}, we obtain:
{ 
\beq\label{B}
\sum_{k=1}^{n} (f_{ij}^{k}-f_{ji}^{k})c_{ml}^{k} =
\sum_{k=1}^{n}\left\lbrace c_{kl}^{i}(f_{jk}^{m}-f_{kj}^{m}) -c_{kl}^{j}(f_{ik}^{m}-f_{ki}^{m}) +c_{mk}^{i}(f_{jk}^{l}-f_{kj}^{l}) - 
c_{mk}^{j}(f_{ik}^{l}-f_{ki}^{l})\right\rbrace.
\eeq
}
Now, let us find the relations associated to the equations~\eqref{eqt1}~-~\eqref{eqt2} of the matched pair of Lie algebras $\G(\A)$ and $\G(\A^{*})$. 
We have $\forall  \, i, j, k: $
\beqs
\left< (-\ad_{\cdot}^{*})(e_i)e_j^{*}, e_k\right>
&=& -\left< \sum_{k=1}^{m}(c_{ik}^{j}-c_{ki}^{j})e_k^{*}, e_k\right>,
\eeqs
providing 
\beq
\displaystyle (-\ad_{\cdot}^{*})(e_i)e_j^{*} = - \sum_{k=1}^{n}(c_{ik}^{j}-c_{ki}^{j})e_k^{*}.
\eeq
Similarly,
\beq
\displaystyle (-\ad_{\circ}^{*})(e_i^{*})e_j = -\sum_{k=1}^{n}(f_{ik}^{j}-f_{ki}^{j})e_k,
\eeq
\beqs
(-\ad_{\circ}^{*})(e_m^{*})[e_i, e_j] 
= \sum_{k=1}^{n}(c_{ij}^{k}-c_{ji}^{k})(-\ad_{\circ}^{*})(e_m^{*})e_k 
 -\sum_{l=1}^{n}\sum_{k=1}^{n}(c_{ij}^{k}-c_{ji}^{k})(f_{ml}^{k}-f_{lm}^{k})e_l.
\eeqs
Then, 
\beq
(-\ad_{\circ}^{*})(e_m^{*})[e_i, e_j] = -\sum_{l=1}^{n}\sum_{k=1}^{n}(c_{ij}^{k}-c_{ji}^{k})(f_{ml}^{k}-f_{lm}^{k})e_l,
\eeq
\beqs
&&-\ad_{\circ}^{*}(\ad_{\cdot}^{*}(e_i)e_m^{*})e_j-[e_i, \ad_{\circ}^{*}(e_m^{*})e_j]+\ad_{\circ}^{*}(\ad_{\cdot}^{*}(e_j)e_m^{*})-
[\ad_{\circ}^{*}(e_m^{*})e_i, e_j] 
\cr&=&
\sum_{l=1}^{n}\sum_{k=1}^{n}\{-(c_{ik}^{m}-c_{ki}^{m})(f_{kl}^{j}-f_{lk}^{j})-(f_{mk}^{j}-f_{km}^{j})(c_{ik}^{l}-c_{ki}^{l}) +
\cr&&
\qquad \qquad (c_{jk}^{m}-c_{kj}^{m})(f_{kl}^{i}-f_{lk}^{i})-(f_{mk}^{i}-f_{km}^{i})(c_{kj}^{l}-c_{jk}^{l})\} e_l. 
\eeqs
Using the fact that $(\G(\A), \G(\A^{*}), \ad_{\cdot}^{*}, \ad_{\circ}^{*})$ is a bimodule of Lie algebras, we have
\beq
\sum_{k=1}^{n}(c_{ij}^{k}-c_{ji}^{k})(f_{ml}^{k}-f_{lm}^{k}) = 
\sum_{k=1}^{n}-(c_{ik}^{m}-c_{ki}^{m})(f_{kl}^{j}-f_{lk}^{j})-(f_{mk}^{j}-f_{km}^{j})(c_{ik}^{l}-c_{ki}^{l})+ \cr 
(c_{jk}^{m}-c_{kj}^{m})(f_{kl}^{i}-f_{lk}^{i})+(f_{mk}^{i}-f_{km}^{i})(c_{jk}^{l}-c_{kj}^{l}),
\eeq
that is, 
{ 
\beqs
\sum_{k=1}^{n}(c_{ij}^{k}-c_{ji}^{k})f_{ml}^{k}+\sum_{k=1}^{n}(c_{ik}^{m}-c_{ki}^{m})f_{kl}^{j}+(c_{ik}^{l}-c_{ki}^{l})f_{mk}^{j}-
(c_{jk}^{m}-c_{kj}^{m})f_{kl}^{i}-(c_{jk}^{l}-c_{kj}^{l})f_{mk}^{i}
= \cr
\sum_{k=1}^{n}(c_{ij}^{k}-c_{ji}^{k})f_{lm}^{k}+\sum_{k=1}^{n}(c_{ik}^{m}-c_{ki}^{m})f_{lk}^{j}+(c_{ik}^{l}-c_{ki}^{l})f_{km}^{j}-
(c_{jk}^{m}-c_{kj}^{m})f_{lk}^{i}-(c_{jk}^{l}-c_{kj}^{l})f_{km}^{i}.
\eeqs
}
 Replacing $l$ (resp. $m$)  by $m$ (resp. $l$)   in the right-hand side of the equality leads to:
 { 
\beq \label{C}
\sum_{k=1}^{n}(c_{ij}^{k}-c_{ji}^{k})f_{ml}^{k}= \sum_{k=1}^{n}\{-(c_{ik}^{m}-c_{ki}^{m})f_{kl}^{j}-(c_{ik}^{l}-c_{ki}^{l})f_{mk}^{j}+
(c_{jk}^{m}-c_{kj}^{m})f_{kl}^{i}+ 
(c_{jk}^{l}-c_{kj}^{l})f_{mk}^{i}\},
\eeq
}
which is identical to the  equation~\eqref{A}.    Besides, 
\beq \label{matched1}
(-\ad_{\cdot}^{*})(e_m)[e_i^{*}, e_j^{*}] =- \sum_{l=1}^{n}\sum_{k=1}^{n}\left\lbrace(f_{ij}^{k}-f_{ji}^{k})(c_{ml}^{k}-c_{lm}^{k})\right\rbrace e_l^{*}
\eeq
\beqs
-\ad_{\cdot}^{*}(\ad_{\circ}^{*}(e_i^{*})e_m)e_j^{*}-[e_i^{*}, \ad_{\cdot}^{*}(e_m)e_j^{*}]+
\ad_{\cdot}^{*}(\ad_{\circ}^{*}(e_j^{*})e_m)e^{*}_m-[\ad_{\cdot}^{*}(e_m)e_i^{*}, e_j^{*}]= \cr
\sum_{l=1}^{n}\sum_{k=1}^{n}\{-(f_{ik}^{m}-f_{ki}^{m})(c_{kl}^{j}-c_{lk}^{j})-(c_{mk}^{j}-c_{km}^{j})(f_{ik}^{l}-f_{ki}^{l}) + 
(f_{jk}^{m}-f_{kj}^{m})(c_{kl}^{i}-c_{lk}^{i})- \\
(c_{mk}^{i}-c_{km}^{i})(f_{kj}^{l}-f_{jk}^{l})\}e_l^{*}.
\eeqs
Then, with $\G(\A) \bowtie^{-\ad_{\cdot}^{*}}_{-\ad_{\circ}^{*}}\G(\A^{*})$ and the relation~\eqref{matched1}, we obtain
\beqs
\sum_{k=1}^{n}(f_{ij}^{k}-f_{ji}^{k})(c_{ml}^{k}-c_{lm}^{k}) = \sum_{k=1}^{n}-(f_{ik}^{m}-f_{ki}^{m})(c_{kl}^{j}-c_{lk}^{j})-
(c_{mk}^{j}-c_{km}^{j})(f_{ik}^{l}-f_{ki}^{l})+\cr
+(f_{jk}^{m}-f_{kj}^{m})(c_{kl}^{i}-c_{lk}^{i})+(c_{mk}^{i}-c_{km}^{i})(f_{jk}^{l}-f_{kj}^{l}),
\eeqs
i.e.,
{ 
\beqs
\sum_{k=1}^{n}(f_{ij}^{k}-f_{ji}^{k})c_{ml}^{k}+\sum_{k=1}^{n}c_{kl}^{j}(f_{ik}^{m}-f_{ki}^{m})+c_{kl}^{i}(f_{jk}^{m}-f_{kj}^{m})-
c_{mk}^{j}(f_{ik}^{l}-f_{ki}^{l}) 
 -c_{mk}^{i}(f_{jk}^{l}-f_{kj}^{l})=\cr 
\sum_{k=1}^{n}(f_{ij}^{k}-f_{ji}^{k})c_{lm}^{k}+\sum_{k=1}^{n}c_{lk}^{j}(f_{ik}^{m}-f_{ki}^{m})+c_{lk}^{i}(f_{jk}^{m}-f_{kj}^{m})-
c_{km}^{j}(f_{ik}^{l}-f_{ki}^{l}) 
 - c_{km}^{i}(f_{jk}^{l}-f_{kj}^{l}),
\eeqs
}
 Replacing $l,$ (resp. $m$),  by $m,$(resp. $l$), in the right-hand side  of the equality leads to
\beq
\sum_{k=1}^{n}(f_{ij}^{k}-f_{ji}^{k})c_{ml}^{k}=\sum_{k=1}^{n}-c_{kl}^{j}(f_{ik}^{m}-f_{ki}^{m})+c_{kl}^{i}(f_{jk}^{m}-f_{kj}^{m})- 
c_{mk}^{j}(f_{ik}^{l}-f_{ki}^{l})+ 
 c_{mk}^{i}(f_{kj}^{l}-f_{jk}^{l}),
\eeq
which  is identical  to the equation~\eqref{B}.
 {$\cqfd$}
 \section*{Aknowledgement}
 This work is partially supported by the Abdus Salam International Centre for Theoretical
 Physics (ICTP, Trieste, Italy) through the Office of External Activities (OEA) - Prj-15. The
 ICMPA is in partnership with the Daniel Iagolnitzer Foundation (DIF), France.

\end{document}